\documentclass{amsart}
\usepackage[latin1]{inputenc}
\newtheorem{thm}{Theorem}
\newtheorem{theorem}[thm]{Theorem}

\newtheorem{lemma}[thm]{Lemma}
\newtheorem{proposition}[thm]{Proposition}
\theoremstyle{definition}
\newtheorem{definition}[thm]{Definition}
\theoremstyle{remark}

\numberwithin{equation}{section}

\newcommand{\To}{\longrightarrow}

\begin{document}

\title{Automatic norm continuity of weak$^\ast$ homeomorphisms}
\author{Antonio Avil\'es}\thanks{The author was supported
by a Marie Curie Intra-European Felloship MCEIF-CT2006-038768
(E.U.) and research projects MTM2005-08379 and S\'{e}neca
00690/PI/04 (Spain)} \address{University of Paris 7, Equipe de
Logique Mathématique, UFR de Mathématiques, 2 place Jussieu, 75251
Paris, France} \email{avileslo@um.es}


\begin{abstract}
We prove that in a certain class $\mathcal E$ of nonseparable
Banach spaces the norm topology of the dual ball is definable in
terms of its weak$^\ast$ topology. Thus, if $X,Y\in\mathcal E$ and
$f:B_{X^\ast}\To B_{Y^\ast}$ is a weak$^\ast$-to-weak$^\ast$
homeomorphism, then $f$ is automatically norm-to-norm
continuous.\\
\end{abstract}

\maketitle

\section{Introduction}

The aim of this note is to prove the automatic continuity in the
norm topology for the weak$^\ast$ homeomorphisms of the dual ball
of certain nonseparable Banach spaces. We start by observing that
such a
property never holds in the separable case.\\

\begin{proposition}\label{casoseparable}
Let $X$ be a separable infinite dimensional Banach space. Then,
there exists a weak$^\ast$-to-weak$^\ast$ homeomorphism
$f:B_{X^\ast}\To B_{X^\ast}$ which is not norm-to-norm continuous.\\
\end{proposition}

Proof: The space $(B_{X^\ast},w^\ast)$ is a metrizable
infinite-dimensional compact convex set, so it is homeomorphic to
the Hilbert cube $[0,1]^\omega$, by Keller's Theorem. It is a
known fact, cf.~\cite[p. 261]{vanMillPrerequisites}, that the
Hilbert cube is countable dense homogenous, which means that if
$A$ and $B$ are countable weak$^\ast$ dense subsets of
$B_{X^\ast}$ then there is a weak$^\ast$ homeomorphism
$f:B_{X^\ast}\To B_{X^\ast}$ such that $f(A)=B$. Thus, it is
enough to find two such subsets $A$ and $B$ with some different
property relative to the norm topology. For example, $A$ can be
taken so that its norm-closure is connected (by choosing it to be
rationally convex) and $B$ with disconnected norm-closure (take
$B'$ a countable weak$^\ast$ dense set which is not norm dense and
then add to $B'$
a point out of its norm-closure).$\qed$\\

We shall introduce two classes $\mathcal E$ and $\mathcal E_0$ of
nonseparable Banach spaces with the properties indicated in the
two
following theorems.\\

\begin{theorem}\label{TheoremClassE}
Let $X$ and $Y$ be spaces in the class $\mathcal E$ and let
$f:B_{X^\ast}\To B_{Y^\ast}$ be a weak$^\ast$-to-weak$^\ast$
homeomorphism. Then $f$ is norm-to-norm continuous.\\
\end{theorem}

\begin{theorem}\label{TheoremClassE0}
Let $X$ and $Y$ be spaces in the class $\mathcal E_0$ and let
$f:B_{X^\ast}\To B_{Y^\ast}$ be a weak$^\ast$-to-weak$^\ast$
homeomorphism. Then, if $(x_{n}^\ast)$ is a sequence in
$B_{X^\ast}$ which
weak$^\ast$ converges to $x^\ast$ and $\|x_n^\ast\|\To \|x^\ast\|$, then $\|f(x_n^\ast)\|\To \|f(x^\ast)\|$.\\
\end{theorem}

The definition of these classes requires some preliminary work and
will be given later but we can already indicate some examples of
spaces which we know that belong in there. The most significant
representatives in $\mathcal E\cap \mathcal E_0$ are the spaces
$c_0(\Gamma)$ and $\ell_p(\Gamma)$ for $1<p<\infty$, $\Gamma$
being an uncountable set of indices. More generally, any
uncountable $c_0$-sum or $\ell_p$-sum ($1<p<\infty$) of separable
spaces belongs to $\mathcal E_0$, and any space from $\mathcal
E_0$ with the dual Kadec-Klee property belongs to $\mathcal E$.
Recall that $X$ has the dual Kadec-Klee property if whenever we
have a weak$^\ast$ convergent sequence $(x_n^\ast)$ in the dual
such that the sequence of norms $(\|x^\ast_n\|)$ converges to the
norm of the limit, then actually the sequence is norm-convergent.
In addition, the class $\mathcal{E}$ is closed under finite
$\ell_1$-sums. On the other hand, one can show that
$\ell_1(\Gamma)\not\in \mathcal E\cup \mathcal
E_0$, indeed:\\

\begin{proposition}
For any infinite set $\Gamma$, there is a
weak$^\ast$-to-weak$^\ast$ homeomorphism
$f:B_{\ell_1(\Gamma)^\ast}\To B_{\ell_1(\Gamma)^\ast}$ which is
not norm-to-norm continuous and does not preserve the limit of norms of weak$^\ast$ convergent sequences.\\
\end{proposition}

Proof: We know that $\ell_1(\Gamma)^\ast = \ell_\infty(\Gamma)$
and hence $B_{\ell_1(\Gamma)^\ast}\approx [-1,1]^\Gamma$ where the
weak$^\ast$ topology is identified with the pointwise topology.
Consider elements $e_1 = (1,0,0,\ldots)$, $e_2=(0,1,0,0,\ldots)$,
etc. in $[-1,1]^\Gamma$. We can easily define coordinatewise a
homeomorphism $f:[-1,1]^\Gamma\To[-1,1]^\Gamma$ taking
$f(\frac{1}{2^n}e_n) = \frac{1}{2}e_n$ and $f(0)=0$. The sequences
$(\frac{1}{2^n}e_n)$ and its image $(\frac{1}{2}e_n)$ weak$^\ast$
converge to 0, but the first is norm convergent and the other is
not,
not even the sequence of norms converges to 0.$\qed$\\

Apart from this extreme example, in fact it happens that these
properties are sensitive to renormings: spaces like
$c_0\oplus_{\ell_1} c_0(\Gamma)$ and $\ell_p\oplus_{\ell_1}
\ell_p(\Gamma)$ fail this automatic norm-continuity property
despite the fact that they are isomorphic to $c_0(\Gamma)$ and
$\ell_p(\Gamma)$, $1<p<\infty$.\\

\begin{proposition}
Let $S$ be an infinite dimensional separable Banach space and $Y$
be a nonseparable Banach space. Let $X=S\oplus_{\ell_1} Y$. Then,
there exists a weak$^\ast$-to-weak$^\ast$ homeomorphism
$f:B_{X^\ast}\To B_{X^\ast}$ which is not norm-to-norm
continuous.\\
\end{proposition}

Proof: Notice that $B_{X^{\ast}} = B_{S^{\ast}}\times B_{Y^\ast}$
and then apply Proposition~\ref{casoseparable}.$\qed$\\

Let us point out where the difficulties appear in proving
Theorem~\ref{TheoremClassE} for say the space $X=c_0(\Gamma)$.
Suppose that $f:B_{X^\ast}\To B_{X^\ast}$ is a weak$^\ast$
homeomorphism. It is a well known fact that the $G_\delta$ points
of $B_{X^\ast}$ are exactly the points of the sphere $S_{X^\ast}$,
hence $f(S_{X^\ast})=S_{X^\ast}$. We know in addition that at the
points of the sphere the norm and weak$^\ast$ topology coincide,
so we conclude that $f$ is norm-continuous at all the points of
the sphere. And this is all that the usual standard
functional-analytic techniques can say to us. In order to get
norm-continuity at all points we shall a need a more powerful tool
coming from topology: Shchepin's spectral theory. This will allow
us to define a certain notion of convergence of a sequence in a
compact space of uncountable weight, that we called
fiber-convergence. Of course, fiber convergent sequences are
respected by homeomorphisms, though not by general continuous
functions. Banach spaces in class $\mathcal E$ are those for which
the fiber convergent sequences in $(B_{X^\ast},w^\ast)$ are
exactly the norm convergent sequences, while the class $\mathcal
E_0$ consists of those spaces in which the fiber convergent
sequences of the dual ball are those weak$^\ast$ convergent
sequences whose sequence of norms converges
to the norm of the limit.\\

\section{Spectral theory}\label{sectionspectral}

In this section, we summarize in a self-contained way what we need
about spectral theory, in the same way as it is exposed in our
joint work with Kalenda~\cite{fiberorders}, which at the same time
is a reformulation of the ideas from \cite{Shchepin1} and
\cite{Shchepin2} in a suitable language. Although this preliminary
material appears already in \cite{fiberorders} with more details, we found it convenient to reproduce it here.\\

Let $K$ be a compact space. We denote by $\mathcal{Q}(K)$ the set
of all Hausdorff quotient spaces of $K$, that is the set of all
Hausdorff compact spaces of the form $K/E$ endowed with the
quotient topology, for $E$ an equivalence relation on $K$. An
element of $\mathcal{Q}(K)$ can be represented either by the
equivalence relation $E$, or by the quotient space $L=K/E$
together
with the canonical projection $p_L:K\To L$.\\

On the set $\mathcal{Q}(K)$ there is a natural order relation. In
terms of equivalence relations $E\leq E'$ if and only if
$E'\subset E$. Equivalently, in terms of the quotient spaces,
$L\leq L'$ if and only if there is a continuous surjection
$q:L'\To L$ such that $q p_{L'} = p_L$. The set $\mathcal{Q}(K)$
endowed with this order relation is a complete lattice, that is,
every subset has a least upper bound or supremum: if $\mathcal{F}$
is a family of equivalence relations of $\mathcal{Q}(K)$, its
least upper bound is the relation given by $x E_0 y$ iff $x E y$
for all $E\in\mathcal{E}$, in other words $E_0 =\sup\mathcal{F} =
\bigcap \mathcal{F}$. It is easy to check that $E_0$ gives a
Hausdorff quotient if each element of $\mathcal{F}$ does.\\

 Let
$\mathcal{Q}_\omega(K)\subset \mathcal{Q}(K)$ be the family of all
quotients of $K$ which have countable weight. Notice that
$\sup\mathcal{A}\in \mathcal{Q}_\omega(K)$ for every countable
subset $\mathcal{A}\subset\mathcal{Q}_\omega(K)$ and also that
$\sup \mathcal{Q}_\omega(K)=K$. A family $\mathcal{S}\subset
\mathcal{Q}_\omega(K)$ is called cofinal if for every $L\in
\mathcal{Q}_\omega(K)$ there exists $L'\in\mathcal{S}$ such that
$L\leq L'$. The family $\mathcal{S}$ is called a
$\sigma$-semilattice if for every countable subset
$\mathcal{A}\subset\mathcal{S}$, the least upper bound of
$\mathcal{A}$ belongs to $\mathcal{S}$.\\

\begin{theorem}[A version of Shchepin's spectral theorem]\label{ourspectraltheorem}
Let $K$ be a compact space of uncountable weight and let
$\mathcal{S}$ and $\mathcal{S'}$ two cofinal $\sigma$-semilattices
in $\mathcal{Q}_\omega(K)$. Then $\mathcal{S}\cap\mathcal{S'}$ is
also a cofinal $\sigma$-semilattice in $\mathcal{Q}_\omega(K)$.\\
\end{theorem}


It is not so obvious to check whether a given $\sigma$-semilattice
is cofinal, so this theorem must be applied together with the
following criterion:\\

\begin{lemma}\label{spectralfactor}
Let $K$ be a compact space of uncountable weight and $\mathcal{S}$
a $\sigma$-semilattice in $\mathcal{Q}_\omega(K)$. Then,
$\mathcal{S}$
is cofinal if and only if $\sup\mathcal{S} = K$.\\
\end{lemma}


The importance of this machinery is that it allows one to study a
compact space of uncountable weight through the study of a cofinal
$\sigma$-semilattice of metrizable quotients and, in particular,
through the natural projections between elements of the
$\sigma$-semilattice. In this way, the study of compact spaces of
uncountable weight is related to the study of continuous
surjections between compact spaces of countable weight. The
following language will be
useful:\\

\begin{definition}
Let $K$ be a compact space of uncountable weight and let
$\mathcal{P}$ be a property. We say that the $\sigma$-typical
surjection of $K$ satisfies property $\mathcal P$ if there exists
a cofinal $\sigma$-semilattice $\mathcal{S}\subset
\mathcal{Q}_\omega(K)$ such that for every $L\leq L'$ elements of
$\mathcal{S}$, the natural projection $p:L'\To L$ satisfies
property
$\mathcal{P}$.\\
\end{definition}

The spectral theorem has the following consequence: In order to
check whether the $\sigma$-typical surjection of $K$ has a certain
property, it is enough to do it on any given cofinal $\sigma$-semilattice, namely:\\

\begin{theorem}
Let $K$ be a compact space of uncountable weight, let
$\mathcal{P}$ be a property, and let $\mathcal{S}$ be a fixed
cofinal $\sigma$-semilattice in $\mathcal{Q}_\omega(K)$. Then the
$\sigma$-typical surjection of $K$ has property $\mathcal{P}$ if
and only if there exists a cofinal $\sigma$-semilattice
$\mathcal{S}'\subset\mathcal{S}$ such that for every $L\leq L'$
elements of $\mathcal{S}'$, the natural projection $p:L'\To L$
satisfies property
$\mathcal{P}$.\\
\end{theorem}

When the compact space we are dealing with is the dual unit ball
$B_{X^\ast}$ of a nonseparable Banach $X$ in the weak$^\ast$
topology, then a cofinal $\sigma$-semilattice in
$\mathcal{Q}_\omega(B_{X^\ast})$ can be obtained from a suitable
family of separable subspaces of $X$.\\

\begin{proposition}\label{Sectionballspectrum}
Let $\mathcal{F}$ be a family of separable subspaces of $X$ such
that
\begin{enumerate}
\item $\overline{span}(\bigcup\mathcal{F}) = X$, and \item if
$\mathcal{F'}\subset\mathcal{F}$ is a countable subfamily, then
$\overline{span}(\bigcup \mathcal{F'})\in \mathcal{F}$.
\end{enumerate}

Then the family $\{B_{Y^\ast} : Y\in\mathcal{F}\}$ is a cofinal
$\sigma$-semilattice in $\mathcal{Q}_\omega(B_{X^\ast})$.
\end{proposition}

Notice that we view $B_{Y^\ast}$ as a quotient of $B_{X^\ast}$ for
$Y\subset X$ through the natural restriction map. The proof of the
proposition is straightforward. In particular, cofinality follows
from the first condition and Lemma~\ref{spectralfactor}.\\

\section{Fiber convergence}

\begin{definition}
Let $\pi:K\To L$ be a continuous surjection, and let $(x_n)$ be a
sequence of elements of $L$ converging to $x\in L$. We say that
$x_n$ is $\pi$-fiber convergent if for every $y\in \pi^{-1}(x)$
there exist elements $y_n\in \pi^{-1}(x_n)$ such that $y_n$
converges to
$y$.\\
\end{definition}

\begin{definition} Let $K$ be a compact space of uncountable weight and let
$(x_n)$ be a sequence in $K$ that converges to $x\in K$. We say
that the sequence $(x_n)$ is fiber convergent if for the
$\sigma$-typical surjection $\pi:L\To L'$, the image of the
sequence in $L'$, $(\pi_{L'}(x_n))_{n<\omega}$, is $\pi$-fiber
convergent.\\
\end{definition}

\begin{definition}
A nonseparable Banach space belongs to the class $\mathcal E$ if
the fiber convergent sequences of $(B_{X^\ast},w^\ast)$ are
exactly the norm
convergent sequences.\\
\end{definition}

\begin{definition}
A nonseparable Banach space belongs to the class $\mathcal E_0$ if
the fiber convergent sequences of $(B_{X^\ast},w^\ast)$ are
exactly those sequences $(x_n^\ast)$ weak$^\ast$ convergent to a
point $x^\ast\in
B_{X^\ast}$ such that $\|x_n^\ast\|\To \|x^\ast\|$.\\
\end{definition}

Notice that Theorems \ref{TheoremClassE} and \ref{TheoremClassE0}
are immediate consequence of the definitions, because the notion
of a fiber-convergent sequence is an intrinsic topological notion
and hence it is preserved under homeomorphisms.\\

\begin{lemma}\label{pifiberosum}
Let $X$ and $Z$ be Banach spaces and $1\leq p \leq\infty$. Set
$Y=X\oplus_{\ell_p} Z$ and let $\pi:B_{Y^\ast}\To B_{X^\ast}$ be
the restriction
map dual to the natural inclusion $X\subset Y$.\\

\begin{itemize}
\item If $p=1$ then every weak$^\ast$ convergent sequence in
$B_{X^\ast}$ is $\pi$-fiber convergent.\\

\item Suppose that $1<p\leq\infty$ and $(x_n^\ast)$ is a sequence
in $B_{X^\ast}$ that weak$^\ast$ converges to $x_0^\ast$. Then
$(x_n^\ast)$ is $\pi$-fiber convergent if and only if the sequence
of norms $(\|x_n^\ast\|)$
converges to $\|x_0^\ast\|$.\\
\end{itemize}
\end{lemma}

Proof: Let $(x_n^\ast)$ be a sequence in $B_{X^\ast}$ that
weak$^\ast$ converges to $x_0^\ast$, and let $y_0^\ast =
x_0^\ast+z_0^\ast\in\pi^{-1}(x_0^\ast)$. If $p=1$, then
$\|x^\ast+z^\ast\| = \max(\|x^\ast\|,\|z^\ast\|)$ for every
$x^\ast\in X^\ast$ and $z^\ast\in Z^\ast$, hence it is enough to
take $y_n^\ast = x_n^\ast+z_0^\ast$ to realize that $(x_n^\ast)$
is actually $\pi$-fiber convergent. If $1<p\leq\infty$, then norms
in the dual are computed as
$$\|x^\ast+z^\ast\| = (\|x^\ast\|^q+\|z^\ast\|^q)^{\frac{1}{q}},\ \ p^{-1}+q^{-1}=1$$

Suppose that the sequence of norms $(\|x_n^\ast\|)$ converges to
$\|x_0^\ast\|$ and let $y_0^\ast = x_0^\ast + z_0^\ast$ be an
arbitrary element of the fiber of $x_0^\ast$. We will find
elements $y_n^\ast\in \pi^{-1}(x_n^\ast)$ such that $y_n^\ast\To
y_0^\ast$. If $z_0^\ast =0$, then we can simply take
$y_n^\ast=x_n^\ast$. Thus, we suppose that $z_0^\ast\neq 0$ and we
define
$$\lambda_n = \max\{\lambda\in [0,1] : \|x_n^\ast + \lambda
z_0^\ast\|\leq 1\}$$ and $y_n^\ast = x_n^\ast + \lambda_n
z_0^\ast$. We have to check that $\lambda_n\To 1$. Suppose on the
contrary that for some subsequence and some $\mu<1$ we have
$\lambda_{n_k} < \mu$. Then
$$(\|x_{n_k}^\ast\|^q+\mu^q\|z_0^\ast\|^q)^{\frac{1}{q}} =
\|x_{n_k}^\ast + \mu z_0^\ast\| > 1$$ so passing to the limit
$$\|x_0^\ast + z_0^\ast\| =(\|x_0^\ast\|^q+\|z_0^\ast\|^q)^{\frac{1}{q}}>(\|x_0^\ast\|^q+\mu^q\|z_0^\ast\|^q)^{\frac{1}{q}} \geq 1$$ which is a contradiction.\\

Conversely, assume now that the sequence of norms $\|x_n^\ast\|$
does not converge to $\|x_0^\ast\|$. Passing to a subsequence, we
can suppose without loss of generality that there is a number
$\mu$ such that $\|x_0^\ast\|<\mu\leq \|x_n^\ast\|$ for every $n$.
Let $\xi\in [0,1]$ be such that $\|x_0^\ast\|^q + \xi^q = 1$, and
let $z_0^\ast$ be any vector of $Z^\ast$ of norm $\xi$. We claim
that there is no sequence $(x_n^\ast + z_n^\ast)\subset
B_{Y^\ast}$ that converges to $x_0^\ast+z_0^\ast$. If it were the
case, then $z_n^\ast\To z_0^\ast$, so $\sup\{\|z_n^\ast\| :
n\in\omega\}\geq \|z_0^\ast\| = \xi$, so $$\sup\{
\|x_n^\ast+z_n^\ast\| : n\in\omega\}\geq \left(\mu^q +
\xi^q\right)^{\frac{1}{q}} > \left(\|x_0^\ast\|^q +
\xi^q\right)^{\frac{1}{q}} = 1$$ which contradicts that $x_n^\ast
+ z_n^\ast\in B_{Y^\ast}$ for every $n$.$\qed$\\

\begin{theorem}\label{mainexamples}
Let $\{X_\alpha : \alpha\in A\}$ be an uncountable family of
separable Banach spaces. \begin{enumerate} \item The $c_0$-sum
$\bigoplus_{c_0}\{X_\alpha : \alpha\in A\}$ belongs to $\mathcal
E_0$, and if $1<p<\infty$ then also $\bigoplus_{\ell_p}\{X_\alpha
: \alpha\in A\}$ belongs to $\mathcal E_0$.
\item All the weak$^\ast$ convergent sequences of the dual ball of $\bigoplus_{\ell_1}\{X_\alpha : \alpha\in A\}$ are fiber convergent.\\
\end{enumerate}
\end{theorem}

Proof: Let $X_0 = \bigoplus_{\alpha\in A} X_\alpha$, where the
type of direct sum is the suitable one in each case. Consider $I$
the family of all countable subsets of $A$ and set $X_i =
\bigoplus_{\alpha\in i} X_\alpha$ for $i\in I$, and
$\mathcal{F}=\{X_i\}_{i\in I}$. This family satisfies conditions
(1) and (2) in Proposition~\ref{Sectionballspectrum}, and
therefore $\mathcal{S} = \{B_{Y^\ast} : Y\in\mathcal{F}\}$ is a
cofinal $\sigma$-semilattice in
$\mathcal{Q}_\omega(B_{X_0^\ast})$. We notice that all the natural
projections between elements of $\mathcal{S}$ correspond to the
dual restriction map of an inclusion of Banach spaces of type
$X\subset X\oplus_{\ell_p} Z$, so that Lemma~\ref{pifiberosum}
indicates which are exactly the $\pi$-fiber convergent sequences
in all those cases. Let us focus on part (1). Let
$(x^\ast_n)\subset B_{X^\ast}$ be a sequence weak$^\ast$
convergent to $x_0^\ast$. For $i\in I$ we denote by
$\pi_i:B_{X_0^\ast}\To B_{X_i^\ast}$ the natural surjection dual
to the inclusion $X_i\subset X_0$. Let $k$ be a countable subset
of $A$ such that $\|x_n^\ast\| = \|\pi_k(x_n^\ast)\|$ for all $n$.
We have then that $\|x_n^\ast\|\To \|x_0^\ast\|$ if and only if
$\|\pi_j(x_n^\ast)\|\To \|\pi_j(x_0^\ast)\|$ for all $j\supset k$,
and by Lemma~\ref{pifiberosum} if and only if $(x_n^\ast)$ is
fiber
convergent.$\qed$\\

Let $KK^\ast$ denote the class of Banach spaces with the dual
Kadec-Klee property. Then, notice that $\mathcal{E}_0\cap \mathcal
E = \mathcal E_0\cap KK^\ast$. Since $c_0(\Gamma)$ and
$\ell_p(\Gamma)$ $(1<p<\infty)$ have $KK^\ast$, we got that these
spaces belong to $\mathcal E$ and satisfy
Theorem~\ref{TheoremClassE}.\\

\begin{lemma}
Let $K=\prod K_i$ be a finite or countable product of compact
spaces of uncountable weight and $(x_n)$ a convergent sequence in
$K$. Then, this sequence is fiber convergent in $K$ if and only if
each of the coordinate sequences $(x_n(i))_{n<\omega}$ is fiber
convergent in
$K_i$.\\
\end{lemma}

Proof: First of all it is straightforward to check that if
$\{\pi_i:X_i\To Y_i\}$ is any family of continuous surjections,
then a sequence $(y_n)\subset \prod Y_i$ is $\prod\pi_i$-fiber
convergent if and only if each coordinate sequence is
$\pi_i$-fiber convergent. Now, going back to the statement of the
lemma, suppose that every coordinate sequence
$(x_n(i))_{n<\omega}$ is fiber convergent. Let $\mathcal{S}_i$ be
a cofinal $\sigma$-semilattice in $\mathcal{Q}_\omega(K_i)$ such
that in all surjections $\pi$ inside $\mathcal{S}_i$, the
projection of the sequence $(x_n(i))$ is $\pi$-fiber convergent.
Let $\mathcal{S}$ be the cofinal $\sigma$-semilattice in
$\mathcal{Q}_\omega(K)$ formed by all quotients of the form $\prod
\pi_i:\prod K_i\To \prod L_i$ where $L_i\in\mathcal{S}_i$. Then,
in all surjections $\pi$ inside $\mathcal{S}$ the projection of
the sequence $(x_n)$ is $\pi$-fiber convergent. Conversely,
suppose that $(x_n)$ is fiber convergent and consider
$\mathcal{S}$ the cofinal $\sigma$-semilattice in
$\mathcal{Q}_\omega(K)$ formed by all quotients which are products
of quotients in each coordinate. There exists a cofinal
$\sigma$-semilattice $\mathcal{T}\subset\mathcal{S}$ such that in
all surjections $\pi$ inside $\mathcal{T}$ the projection of the
sequence $(x_n)$ is $\pi$-fiber convergent. For every $i$, we
consider $\mathcal{T}_i$ to be the set of all quotients $L$ of
$K_i$ such that there is some quotient $\prod_j L_j$ with $L_i=L$.
Then $\mathcal{T}_i$ is a cofinal $\sigma$-semilattice of
$\mathcal{Q}_\omega(K_i)$ and, by the observation at the beginning
of this proof, in all surjections $\pi$ inside $\mathcal{T}_i$ the
projection of $(x_n(i))$ is $\pi$-fiber convergent.$\qed$\\

\begin{proposition}
Let $\{X_n : n<\omega\}$ be a countably infinite family of
nonseparable Banach spaces. Then $\bigoplus_{\ell_1}\{X_n :
n<\omega\}$ belongs neither to $\mathcal E_0$ nor to $\mathcal E$.\\
\end{proposition}

Proof: Notice that $B_{X^\ast} = \prod_{n<\omega}B_{X_n^\ast}$.
For $x_n$ an element of the sphere of $X_n^\ast$, the sequence
$(x_0,0,0,\ldots)$, $(0,x_1,0,\ldots)$, $\ldots$ is fiber
convergent to 0 but not norm
convergent.$\qed$\\

\begin{proposition}\label{el1sum}
If $X,Y\in \mathcal E$, then $X\oplus_{\ell_1} Y\in\mathcal
E$.\\
\end{proposition}

Proof: Let $Z =X\oplus_{\ell_1} Y$. Notice that $B_{Z^\ast} =
B_{X^\ast}\times B_{Y^\ast}$. A sequence in this product is
fiber-convergent if and only if both coordinates are
fiber-convergent and the same happens for
norm-convergence.$\qed$\\

We can provide a couple of extra examples. In both cases $1<p<\infty$ and $c_0$ can be substituted for $\ell_p$:\\

\begin{itemize}

\item $\ell_p(\Gamma)\oplus_{\ell_1} \ell_p(\Gamma)\in \mathcal
E\setminus \mathcal E_0$. This space belongs to $\mathcal E$ by
Proposition~\ref{el1sum}, but it is not in $\mathcal{E}_0$ because
it is in $\mathcal{E}$ but it fails the dual Kadec-Klee property:
If $(x_n^\ast)$ is a sequence inside the unit sphere of
$\ell_p(\Gamma)^\ast$ which weak$^\ast$ converges to 0, then the
sequence $(x_0^\ast,x_n^\ast)$ shows that $KK^\ast$ does not
hold.\\

\item An uncountable $\ell_p$-sum of copies of
$\ell_p(\omega)\oplus_{\ell_1}\ell_p(\omega)$ belongs to $\mathcal
E_0\setminus \mathcal E$. This space belongs to $\mathcal{E}_0$ by
Theorem~\ref{mainexamples} but again it is not in $\mathcal{E}$
since it is in $\mathcal{E}_0$ but it fails $KK^\ast$ for similar
reasons as in the previous case.

\end{itemize}

\end{document}